\documentclass[11pt,a4paper]{amsart}
\usepackage[latin1]{inputenc}

\newtheorem{theo}{Theorem}[section]
\newtheorem{prop}[theo]{Proposition}
\newtheorem{cor}[theo]{Corollary}

\theoremstyle{definition}
\newtheorem{rque}[theo]{Remark}

\usepackage{amssymb}
\usepackage{amsfonts}
\usepackage{amsmath}
\usepackage{tabularx}

\DeclareMathOperator{\F}{\mathbb{F}}
\DeclareMathOperator{\Q}{\mathbb{Q}}

\DeclareMathOperator{\Z}{\mathbb{Z}}

\DeclareMathOperator{\Cl}{Cl} 
\DeclareMathOperator{\Frob}{Frob}
\DeclareMathOperator{\G}{\mathcal{G}} 
\DeclareMathOperator{\GL}{GL}
\DeclareMathOperator{\SL}{SL}

\DeclareMathOperator{\Ind}{Ind}
\DeclareMathOperator{\N}{N} 
\DeclareMathOperator{\Nm}{\N_{F\!/\!\Q}}
\DeclareMathOperator{\PGL}{PGL}
\DeclareMathOperator{\tr}{tr}
\DeclareMathOperator{\rH}{\mathrm{H}}

\newcommand{\cf}{{\it cf }}

\begin{document}

\title[images of  Galois representations attached  to Hilbert modular
forms]{Explicit determination of  images of  Galois representations 
attached  to Hilbert modular forms}

\author{Luis Dieulefait and Mladen Dimitrov}
\date{\today}

\begin{abstract} In \cite{dim-hmv} the second author proved that the 
image of the Galois representation mod $\lambda$  attached to a 
Hilbert modular newform  is ``large" for 
all but finitely many primes $\lambda$, if the form is not a theta
series.  In this brief note, we give an explicit
bound for this exceptional finite set of primes and determine the images in three different
examples. Our examples are of Hilbert newforms on real quadratic fields, of
parallel or non-parallel weight and of different levels.
\end{abstract}

\maketitle

\vspace{1cm}
\section{Introduction}

Serre  \cite{serre2} and Ribet \cite{ribet3,ribet6}   proved a large
image result for the  compatible family of 
Galois representations attached to a classical  modular newform $f$, provided 
that $f$ is not a theta series. The result is the following: let
$E$ be the number field generated by the eigenvalues of $f$ and for
every prime $\lambda$  in $E$, let $\rho_\lambda$ be the continuous 
$\lambda$-adic Galois representation  attached to $f$ constructed by Deligne; 
then for all but finitely many primes  $\lambda$ the image of $\rho_\lambda$
 is large, namely it contains  $\SL_2(\mathbb{Z}_\ell)$, with $\ell$ the rational 
prime below $\lambda$.
  In fact, the exact value of the image for all but finitely many primes  
$\lambda$ can be given 
if the inner-twists  of $f$ (if any) are taken into account.
 The finite set of primes where the above result fails is usually called the 
   ``exceptional set" of the family $\{ \rho_\lambda \}$. 
    Explicit versions of  the result of Ribet are known and the exceptional sets 
  have been explicitly computed
  for several examples of classical modular forms \cite{Di,DV}. 
  We will now recall a result of \cite{dim-hmv} that is the natural generalization
  of Ribet's result to the case of Hilbert modular forms.

\medskip
Let $F$ be a totally real number field of degree $d$.
Let $k=(k_1,..,k_d)$ be an arithmetic weight ({\it i.e.} $k_i\geq 2$ are of the same parity) 
and put $k_0=\max\{k_i|1\leq i\leq d\}$. Let $N$ be an integral ideal of $F$ 
and $\psi$   a Hecke character of $F$  of conductor dividing $N$ and  infinity 
type $2-k_0$.  We consider a Hilbert modular newform $f$ over $F$ of weight $k$, 
level $N$ and central character $\psi$. By a theorem of 
Shimura, the Fourier coefficients $c(f,\pi)$ of $f$ ($\pi$ is a prime of $F$)
generate a number field $E$. 
The  absolute Galois  group of a field $L$ is denoted  $\G_L$.

\medskip
By the work of Ohta, Carayol, Blasius-Rogawski and  Taylor 
\cite{Ta}, for every prime  $\lambda$ of $E$ one  can associate
to $f$ an   absolutely irreducible,   totally odd $\lambda$-adic representation 
$\rho_{\lambda}:\G_F\rightarrow \GL_2(E_\lambda)$,   unramified outside $N l$.
By a theorem of Carayol \cite{carayol}  the restriction of $\rho_{\lambda}$
to the decomposition group at a prime  $\pi$ of $F$ not dividing $\ell$ is determined
by the local Langlands correspondence for $\GL_2$. In particular, for 
all  primes  $\pi$ of $F$ not dividing $N l$ we have :
$$\tr(\rho_\lambda(\Frob_\pi))=c(f,\pi) \qquad
\det(\rho_\lambda(\Frob_\pi))=\psi(\pi)\Nm(\pi),$$
where  $\Frob_\pi$ denotes a geometric Frobenius at $\pi$.
 
\medskip
Denote by $\F_\lambda$ the residue field of $E_\lambda$.
By taking a Galois stable lattice, we define 
$\overline{\rho}_{\lambda}=\rho_{\lambda}
\mod{\lambda}:\G_F \rightarrow \GL_2(\F_\lambda)$,  
whose  semi-simplification is independent of the particular choice of
a lattice.

By \cite[Prop.3.8]{dim-hmv} if $f$ is not a theta series, then 
 for all primes $\lambda$ outside a finite set of ``exceptional primes"
 the image of  $\overline{\rho}_{\lambda}$ contains $\SL_2(\F_{\ell})$.

\medskip
In this article, we  determine explicitly   the images and  the finite exceptional sets
for three examples of  Hilbert modular newforms on real quadratic number fields,
of different weights and levels. Our results can be summarized in the
following table : 

\bigskip
\begin{center}\begin{tabularx}{\linewidth}{|l|c|c|c|c|c|}
\cline{1-6} 
  $f$ constructed by  &   $F$  & $E$  &   k   &
N   &  exceptional $\lambda$ divide  \\ 
\cline{1-6}
 Consani-Scholten  & $\Q(\sqrt{5})$  &   $\Q(\sqrt{5})$  &  
$(2,4)$   & $2\cdot 3\cdot 5 $   & $2\cdot 3\cdot 5$  (Thm.\ref{cs})\\ 
\cline{1-6} 
  Demb{\'e}l{\'e}  &   $\Q(\sqrt{5})$   &   $\Q(\sqrt{6})$  &  
$(2,2)$   & $5\cdot(8+3\sqrt{5}) $ 
& $2\cdot 3\cdot 5\cdot 19$  (Thm.\ref{dembele}) \\ 
\cline{1-6} 
  Okada  &   $\Q(\sqrt{257})$   &   $\Q(\sqrt{13})$  &  
$(2,2)$   & $1$   
& $2\cdot 3\cdot 257 $  (Thm.\ref{okada})\\ 
\cline{1-6} 
\end{tabularx}
\end{center}

\bigskip
Section \ref{modularity} contains a potential application to the
modularity of a quintic threefold.

\section{Large image results for Hilbert modular forms.}

In this section we recall some of the  results of
\cite[\S3]{dim-hmv}. 

\medskip
Let $v$ be a prime of $F$ above $\ell$ of residual degree $h$. 
Put $\overline{\rho}=\overline{\rho}_{\lambda}|_{D_{v}}$. 
The semisimplification $\overline{\rho}^{\mathrm{ss}}$ of 
$\overline{\rho}$ is tamely ramified and the image by
$\overline{\rho}^{\mathrm{ss}}$ of the tame inertia $I_{v}^t$
is cyclic. Let $a$ be a generator of 
$T=\overline{\rho}^{\mathrm{ss}}(I_{v}^t)$. Let $n$ be the 
image of a Frobenuis element by $\overline{\rho}^{\mathrm{ss}}$. 
Then $n$ belongs to the normalizer $N$ of 
$T$ in  $\GL_2(\F_\lambda)$ and  $nan^{-1}=a^{\ell^h}$. 
Since $N/T$ is of order $2$, we have

$\bullet$ either $n\in T$ and $a=a^{\ell^h}$. In this case 
there exist integers $(p_i)_{0\leq i\leq 2h-1}$ such that 
$\overline{\rho}^{\mathrm{ss}}\!:\! I_{v}^t\!\rightarrow
\GL_2(\F_\lambda)$ factors through the natural map 
$I_{v}^t\rightarrow \F_{\ell^h}^\times$ followed by 
$$F_{\ell^h}^\times \rightarrow \F_{\ell^h}^\times\times
\F_{\ell^h}^\times\enspace,\enspace
x\mapsto (x^{p_0+p_1\ell+...+p_{h-1}\ell^{h-1}},x^{p_h+p_{h+1}\ell+...+p_{2h-1}\ell^{h-1}}),$$

$\bullet$ either $n\notin T$ and $a=a^{\ell^{2h}}$. In this case 
there exist integers $(p_i)_{0\leq i\leq 2h-1}$ such that 
$\overline{\rho}^{\mathrm{ss}}\!:\! I_{v}^t\!\rightarrow
\GL_2(\F_\lambda)$ factors through the natural map 
$I_{v}^t\rightarrow \F_{\ell^{2h}}^\times$ followed by 
$$F_{\ell^{2h}}^\times \rightarrow \F_{\ell^{2h}}^\times
\times\F_{\ell^{2h}}^\times \enspace,\enspace
x\mapsto (x^{p_0+p_1\ell+...+p_{2h-1}\ell^{2h-1}},x^{\ell^h(p_0+p_1\ell+...+p_{2h-1}\ell^{2h-1})}).$$
\begin{prop}{\rm\cite[Cor.2.13]{dim-hmv}} \label{weights}
Assume  that  $\ell>k_0$ is unramified in $F$ and  does not divide $N$. 
Then $\overline{\rho}_{\lambda}$ is crystalline at $\ell$, the 
  multisets $\bigcup_{v|\ell} \{p_i| 0 \leq i\leq 2h-1\}$ and 
$\{\frac{k_0-k_i}{2}, \frac{k_0+k_i-2}{2}| 1\leq i\leq d\}$ are equal
and $p_{i+h}=p_i$, for all $0\leq i\leq h-1$.
\end{prop}

Using this proposition, one can show : 

\begin{prop}\label{image}
Assume  that  $\ell\!>\!k_0$ is unramified in $F$ and  does not divide $N$. 

{\rm(i)} {\rm\cite[Prop.3.1(i)]{dim-hmv}}
For all but finitely many primes $\lambda$,
 $\overline{\rho}_{\lambda}$ is absolutely  irreducible.

\smallskip
{\rm(ii)} {\rm\cite[\S3.2]{dim-hmv}} If $d(\ell\!-\!1) >  5\sum_{i=1}^d(k_i\!-\!1)$
then the image of $\overline{\rho}_{\lambda}$ in $\PGL_2(\F_\lambda)$ is
not isomorphic  to one of the groups $A_4$, $S_4$ or $A_5$. 

\smallskip
{\rm(iii)} {\rm\cite[Lemma 3.4]{dim-hmv}}
 Assume that  $\ell\neq 2k_i-1$ for all $1\leq i \leq d$.
Assume that  the image   of $\overline{\rho}_{\lambda}$ in
$\PGL_2(\F_\lambda)$ is dihedral, that is  $\overline{\rho}_{\lambda}\cong
\overline{\rho}_{\lambda}\otimes \varepsilon_{K/F}$ where
$\varepsilon_{K/F}$ denotes the character of some quadratic extension
$K/F$. Then $K/F$ is unramified outside $N$. 
\end{prop}

\begin{cor} {\rm\cite[Prop.3.8]{dim-hmv}}
 If  $f$ is not a theta series then for 
all but finitely many primes $\lambda$ the image of
 $\overline{\rho}_{\lambda}$ contains (a conjugate of) 
$\SL_2(\F_{\ell})$. 
\end{cor}

The proof uses Dickson's classification theorem of the subgroups of 
$\GL_2(\F_\lambda)$. According to this theorem, the image in 
 $\PGL_2(\F_\lambda)$ of an irreducible subgroup of
 $\GL_2(\F_\lambda)$ that does not contain (a conjugate of) 
$\SL_2(\F_{\ell})$ is  isomorphic either to a dihedral group, either to 
one of the groups $A_4$, $S_4$ or $A_5$.

\section{On an example of Consani-Scholten.}\label{ex1}

In \cite{CoSc} Consani and Scholten construct a Hilbert modular newform $f$
on $F=\Q(\sqrt{5})$ of weight $(2,4)$, conductor $30$ and  trivial
central character.
The  Fourier coefficients $c(f,\pi)$ of $f$ belong to $E=\Q(\sqrt{5})$
and  are computed for $7\leq \Nm(\pi)\leq 97$. By Prop.\ref{weights},
for all $\ell\geq 7$,  $\overline{\rho}_{\lambda}$  is crystalline and
the  semi-simplifications of its restrictions to the inertia groups at
$\ell$ satisfy :

$\bullet$ for  $\ell=vv'$  split in $F$  :
$$\overline{\rho}_{\lambda|}{}_{I_v}^{\mathrm{s.s.}}\simeq
 1 \oplus \omega^3 \text{ and }
\overline{\rho}_{\lambda|I_{v'}}^{\mathrm{s.s.}}\simeq
 \omega \oplus \omega^2 \text{, or }
\overline{\rho}_{\lambda|}{}_{I_v}^{\mathrm{s.s.}}\simeq
 \omega_2^3 \oplus \omega_2^{3\ell} \text{ and }
\overline{\rho}_{\lambda|}{}_{I_{v'}}^{\mathrm{s.s.}}\simeq
 \omega_2^{1+2\ell} \oplus \omega_2^{2+\ell} .$$

$\bullet$ for $\ell$ inert  in $F$  :
$$\overline{\rho}_{\lambda|}{}_{I_{\ell}}^{\mathrm{ s.s.}} \simeq
 \omega_2 \oplus \omega_2^{2+3\ell} \text{ or }
 \omega_2^{1+3\ell} \oplus \omega_2^2 \text{ or }
\omega_4^{a+3\ell+(3-a)\ell^2} \oplus \omega_4^{3-a+a\ell^2+3\ell^3},
\text{with } a=1 \text{ or } 2.$$

Here $\omega$ is the cyclotomic character modulo $\ell$ and $\omega_2$ (resp. $\omega_4$)
denotes a fundamental character of level $2$ (resp. $4$).

\medskip
The aim of this section is to  establish  the following

\begin{theo} \label{cs} Let  $\ell=\Nm(\lambda)\geq 7$ be a prime.
Then the image of 
 $\overline{\rho}_{\lambda}$ is isomorphic to
$\begin{cases}
 \{\gamma\in\GL_2(\F_{\ell}),  \det(\gamma)\in \F_{\ell}^{\times 3}\}\text{, if }
\ell\equiv \pm 1\pmod{5}, \\
 \{\gamma\in \GL_2(\F_{\ell^2}),  \det(\gamma)\in \F_{\ell}^{\times 3}\}\text{, if }
\ell\equiv \pm 2\pmod{5}.\end{cases}$
\end{theo}

\noindent {\bf Proof :} The theorem would follow from
 \cite[Prop.3.9]{dim-hmv} once  
we establish that the image of $\overline{\rho}_{\lambda}$ is
not reducible, nor of  dihedral type, nor of $A_4$, $S_4$ or  $A_5$ type, 
nor isomorphic to  
$\{\gamma\in \F_{\ell^2}^\times\GL_2(\F_{\ell})\text{ , } \det(\gamma)\in \F_{\ell}^{\times 3}\}$
(in the last case we say that $\overline{\rho}_{\lambda}$ has  inner twists).

$\relbar$ Denote by $\epsilon=\frac{1+\sqrt{5}}{2}$ the fundamental unit of $F$.
Then $\epsilon^{120}\equiv 1\pmod{30}$.  By
\cite[Prop.3.1(ii)]{dim-hmv} we obtain that $\overline{\rho}_{\lambda}$ is
absolutely irreducible unless $\epsilon^{240}-1\in \lambda$ that is
(recall that $\ell \geq 7$)
\begin{equation}\label{list1}
 \ell\in\{7,11,23,31,61,241,599,1553,2161,20641\}
\end{equation}
For odd representations irreducibility is equivalent to absolute irreducibility.
 Since we  are only interested in showing the irreducibility for the finite set
 of primes in (\ref{list1}) it is enough to prove that for each $\lambda$ in $E$
 dividing such a prime $\ell$  there exists a prime $\pi$ in $F$ relatively prime to
$30\ell$ such that the characteristic polynomial of $\bar{\rho}_{\lambda}
(\Frob_{\pi})$ is irreducible over  $\mathbb{F}_\lambda$.
 We have used the characteristic polynomials for the primes $\pi$
 in $F$ dividing one of the following  primes $p$:
\begin{equation}\label{list2}
 11, 19, 29, 31, 41, 59 
\end{equation}
 Observe that all these primes split in $F$. For each of these primes $\pi$ the
 characteristic polynomial  $x^2- c(f,\pi) x + p^3$ has been 
 computed in \cite{CoSc}. For every prime  $\ell \geq 7$
  in (\ref{list1}) and for every  $\lambda$ in $E$ dividing $\ell$
   we checked that some of these characteristic polynomials are irreducible
  modulo $\lambda$, and this proves that for every prime $\ell \geq 7$ the
  residual representations are absolutely irreducible.

\smallskip
$\relbar$ Assume that the image of $\overline{\rho}_{\lambda}$ is of dihedral type
and $\ell\geq 7$.  By the local behavior of $f$ at $2$ and $3$ (\cf \cite{carayol}) 
and by  Prop.\ref{image}(iii)  there exists a   quadratic  extension $K/F$ of discriminant  dividing
$5$  such that  for every prime $\pi$ of $F$ inert in $K/F$ and
relatively prime to $30\ell$, $\tr(\overline{\rho}_\lambda (\Frob_\pi))=0$.
Equivalently, for every such prime $\pi$ we
 have  $c(f,\pi)\in \lambda$.

Thus, the algorithm for bounding the set of dihedral primes runs as follows: for
each quadratic extension $K$ of $F$ unramified outside $5$ find a couple of primes
$\pi$ in $F$ verifying  $\pi$ inert in $K/F$ and $c(f,\pi) \neq 0$. 
The only possibly
dihedral primes are those $\lambda$ dividing these $c(f,\pi)$. Using again the
eigenvalues $c(f,\pi)$ for the primes $p$ listed in (\ref{list2}) we have checked with this
method that there is no dihedral prime $\ell \geq 7$.

Let us describe in more detail the computations performed : the following list contains all quadratic extensions of
$F = \Q(\sqrt{5})$ unramified outside $\sqrt{5}$ (here we use the fact that the multiplicative group of units of $F$ is generated by
$-1$ and $\epsilon=\frac{1+\sqrt{5}}{2}$) :
$$ K_1 = F (\sqrt[4]{5}); \; K_2 = F\Bigl( \sqrt{  \epsilon \sqrt{5}} \Bigr); \;
K_3 = F\Bigl(\sqrt{-\sqrt{5}}\Bigr); \; K_4 = F\Bigl(\sqrt{- \epsilon \sqrt{5}}\Bigr)$$
\begin{rque} It is easy to see that $K_4$ is the cyclotomic 
field of $5$-th roots of unity.
\end{rque}
The primes $\pi$ in $F$ dividing $p = 29 \text{ or } 41$
are inert in $K_1 / F$ and in $K_3 / F$. 
The primes $\pi$ in $F$ dividing $p = 11 \text{ or } 29$ (resp.
$p = 29 \text{ or } 59$) are inert in $K_2 / F$ (resp. in $K_4 / F$). 
All these primes $\pi$ verify $c(f, \pi) \neq 0$.
In each case, we search for all primes $\lambda$ dividing
$N_p := \N_{E/\Q} c(f, \pi)$ for both primes $p$ in the above pairs.
By computing the prime factors of each of the  $N_p$ :
$$ N_{11} : \{ 2,11,19 \} \enspace\enspace N_{29} : \{ 2,5,29 \} \enspace\enspace
 N_{41} : \{ 2,41,379 \} \enspace\enspace  N_{59} : \{ 2,5,59,71 \} $$
we conclude that there are no dihedral primes $\ell\geq 7$.


$\relbar$ By  Prop.\ref{image}(ii) the  image of $\overline{\rho}_{\lambda}$ cannot be
of $A_4$, $S_4$ or $A_5$ type  for $\ell\geq 11$ and a more careful
study  shows that the argument remains  valid for $\ell=7$.

$\relbar$ Finally, assume that $\overline{\rho}_{\lambda}$ has inner
twists, that is the  image of $\overline{\rho}_{\lambda}$ is isomorphic
to $\{\gamma\in \F_{\ell^2}^\times\GL_2(\F_{\ell})\text{ , } \det(\gamma)\in \F_{\ell}^{\times 3}\}$.
In this case $\ell$ should be  inert in $E$ and  the squares of the traces of the
elements of the image of $\overline{\rho}_{\lambda}$ should belong to $\F_{\ell}$
(\cf \cite[\S3.4]{dim-hmv}). In particular $c(f,\pi)^2\in \F_{\ell}$
for all the primes $\pi$ in $F$  not dividing $30\ell$. Taking $\pi=(5+\epsilon)$ we
find a contradiction. \hfill $\square$

\section{Towards the modularity of a quintic threefold.}\label{modularity}

Here we describe a potential application.  Consani and  Scholten
study  the middle degree {\'e}tale cohomology of a quintic threefold
$\widetilde{X}$
(a proper and smooth $\Z[\frac{1}{30}]$-scheme with Hodge numbers
 $h^{3,0}=h^{2,1}=1$, $h^{2,0}=h^{1,0}=0$ and $h^{1,1}=141$). They
show that the $\G_{\Q}$-representation $\rH^3(\widetilde{X}_{\overline{\Q}},\Q_{\ell})$
is induced from a two dimensional representation $\sigma_\lambda$ of $\G_F$, 
and they conjecture that $\sigma_\lambda$ is isomorphic to the $\rho_{\lambda}$
from section \ref{ex1}.

\begin{prop} \label{modularite} Assume $\ell\geq 7$ and  $\overline{\rho}_{\lambda}=
\sigma_\lambda  \mod{\lambda}$. Then  $\sigma_\lambda$ is modular, and 
in particular the $L$-function associated to 
$\rH^3(\widetilde{X}_{\overline{\Q}},\Q_{\ell})$ has  an analytic continuation 
to the whole complex plane and a functional equation.
\end{prop}

\noindent {\bf Proof : } By Thm.\ref{cs} and  \cite[Prop.3.13]{dim-hmv} the ``large''
image condition $\mathrm{{\bf (LI}_{Ind \overline{\rho}_{\lambda}}{\bf )}}$ 
of  \cite[\S0.1]{dim-hmv}  is fulfilled for all  $\ell\geq 7$. Then  
the proposition follows from  \cite[Thm.0.2]{dim-ihara}.
\hfill $\square$

\begin{rque}
{\rm (i)} In order to prove that $\overline{\rho}_{\lambda}= \sigma_\lambda$ mod $\lambda$ it is
enough, by the Cebotarev Density Theorem, to check that for ``sufficiently many''
primes $\pi$ of $F$ we have $\tr(\sigma_\lambda(\Frob_\pi)) \equiv c(f,\pi) \pmod{\lambda}$.
Unfortunately, a relatively small number of these values have been computed.

{\rm (ii)} If  $\ell$ is inert in $E$, the image of 
$\Ind_F^{\Q} \overline{\rho}_\lambda: \G_F \rightarrow 
(\GL_2\times \GL_2)(\F_{\ell^2})$ is isomorphic to
$$\{(\gamma,\gamma^\sigma)\enspace|\enspace\gamma\in \GL_2(\F_{\ell^2}),\det(\gamma)\in \F_{\ell}^{\times 3}\}$$
where $\sigma$ denotes the non-trivial automorphism of $\F_{\ell^2}$. This is 
a typical example of a representation with non-maximal image satisfying 
$\mathrm{{\bf (LI}_{Ind \overline{\rho}_{\lambda}}{\bf )}}$.

{\rm (iii)} The results of Taylor\cite{Ta3}, if extended to the
non-parallel weights $k$, would give the meromorphic   continuation 
to the whole complex plane and the functional equation.
\end{rque}

\section{On examples of Demb{\'e}l{\'e}.}

In \cite[p.60]{Dem} Demb{\'e}l{\'e}  constructs a  Hilbert modular newform $f=f_4$
on $F=\Q(\sqrt{5})$ of parallel weight $(2,2)$, level $N=\sqrt{5}^2(8+3\sqrt{5})$, and
Fourier coefficients in  $E=\Q(\sqrt{6})$.

By Prop.\ref{weights} for all $\ell$ not dividing $5\cdot 19$,   $\overline{\rho}_{\lambda}$
is crystalline at $\ell$ and 
$$\overline{\rho}_{\lambda|}{}_{I_{\ell}}^{\mathrm{s.s.}}\simeq
 1 \oplus \omega \text{  or }  \omega_2 \oplus \omega_2^{l} \text{  or }
\omega_4^{1+l} \oplus \omega_4^{\ell^2+\ell^3}.$$

\begin{theo} \label{dembele} Assume  $\ell\geq 7$ and $\ell\neq
 19$. Then the image of   $\overline{\rho}_{\lambda}$
is isomorphic to $\begin{cases}  \GL_2(\F_{\ell}) \text{, if } \ell\equiv \pm 1,\pm 5\pmod{24}, \\
 \{\gamma\in \GL_2(\F_{\ell^2}),\det(\gamma)\in\F_{\ell}^{\times }\}\text{, otherwise}.\end{cases}$
\end{theo}

\noindent {\bf Proof : } As in Thm.\ref{cs},  we have to discard  the reducible,  
dihedral, $A_4$, $S_4$, $A_5$ and inner twist cases.

$\relbar$ Assume that  $\overline{\rho}_{\lambda}$ is reducible :
$\overline{\rho}_{\lambda}^{\mathrm{s.s.}}=
\varphi\oplus\omega\varphi^{-1}$, where $\varphi$ is a mod $\lambda$ crystalline character of
$\G_F$. Using the compatibility between local and global Langlands correspondence
for $\GL_2$ (\cf \cite{carayol}) and the fact that $f$ has trivial central character,
we deduce that the prime-to-$\ell$ part of the conductor  of  $\varphi$ divides   $\sqrt{5}$.
It is easy to check that  $\epsilon^{4}\equiv 1\pmod{\sqrt{5}}$ and that 
$\lambda$ does not divide $\epsilon^{4}-1$. By the proof of  \cite[Prop.3.1]{dim-hmv} 
 it follows  that $\varphi$ is unramified at $\ell$,
and therefore corresponds (via global class field theory) to a character of the
ray class group $\Cl_{F,\sqrt{5}}^+$. Now we use that  $\Cl_{F,\sqrt{5}}^+$ is of 
order $2$, because the narrow class group  $\Cl_{F}^+$ is trivial and  $\epsilon^2$
generates an index $2$ subgroup  of $(\Z[\epsilon]/(\sqrt{5}))^\times\cong 
(\F_5[X]/X^2)^\times$.  Hence $\varphi$ is a trivial or quadratic character.  
By evaluating
$\overline{\rho}_{\lambda}^{\mathrm{s.s.}}=\varphi\oplus\omega\varphi^{-1}$ at 
$\Frob_7$ we find a contradiction (the conclusion does not apply to the prime 
$\lambda=7$, however in the next few lines we give a second proof of the residual 
irreducibility that also covers this prime).

Once we know that $\varphi$ is unramified at $\ell$ and its conductor divides $\sqrt{5}$,
an alternative method
consists (without determining the order of $\varphi$) of taking primes that are trivial
modulo $\sqrt{5}$, more precisely :
 find primes $\pi$ in $F$ having a generator $\alpha$ (in $F$ all ideals are principal)
 such that  $\alpha \equiv 1 \pmod{\sqrt{5}}$.
 It follows from class field theory that for all such primes $\pi$, the element $\Frob_\pi$ is totally
 split in the ray class field of $F$ of conductor $\sqrt{5}$, hence
 $\varphi$ is trivial on  $\Frob_\pi$.   This property is fulfilled by the two primes of $F$ 
dividing $p=11$ and $31$. In fact  $ 11 = (-4+\sqrt{5}) (-4-\sqrt{5})$ and $31 = (6+ \sqrt{5}) (6-
 \sqrt{5})$.   Therefore, equating traces in
$\overline{\rho}_{\lambda}^{\mathrm{s.s.}}= \varphi\oplus\omega\varphi^{-1}$
at these $\Frob_\pi$ we obtain :
$$ c(f,-4+\sqrt{5} ) = -2 \sqrt{6} \equiv 1 + 11 \pmod{\lambda}
\text{, } \enspace c(f,6+ \sqrt{5}) = 2 \equiv 1 + 31 \pmod{\lambda}$$
From this, we conclude again that the residual representation is irreducible for 
 $\ell > 5$, $\ell \neq 19$.

$\relbar$  Assume that   $\overline{\rho}_{\lambda}$ is of dihedral type. As in the
previous example, by the local behavior of  $f$ at $ 8 + 3 \sqrt{5}$
and by  Prop.\ref{image}(iii)  we have to consider only  quadratic extensions of $F$ unramified
outside $\sqrt{5}$, that is  the same extensions $K_1,...,K_4$ that we considered in
section \ref{ex1}.

For each $1\leq i\leq 4$, we  take again two rational primes $p$ among $11$, $29$, $41$, $59$
decomposing as $\pi \pi'$ in $F$ with both $\pi$ and $\pi'$ inert in $K_i /F$.
As in the previous example, let us list the prime factors of the norms $N_p = \N_{E/\Q} c(f, \pi) $
for each of these $p$  and a suitably chosen $\pi \mid p$ :
$$ N_{11} : \{ 2,3 \} \qquad N_{29} : \{ 3,5 \}  \qquad 
 N_{41} : \{ 2,3 \} \qquad  N_{59} : \{ 3,5 \} $$

Since no prime $\ell \geq 7$ divides any of these norms, we  conclude that
no prime $\ell\geq 7$, $\ell \neq 19$, is dihedral.

$\relbar$ By Prop.\ref{image}(ii) the image of $\overline{\rho}_{\lambda}$ cannot be
of $A_4$, $S_4$ or $A_5$ type  for any $\ell\geq 7$.

$\relbar$  As in the proof of Thm.\ref{cs} we see 
that $\overline{\rho}_{\lambda}$ cannot have inner twists. \hfill $\square$

\begin{rque} 
We applied this method to bound explicitly  the set of exceptional primes 
for several other examples of Hilbert modular newforms from Demb{\'e}l{\'e}'s tables 
and it always worked well. 
\end{rque}

\section{On an example of Okada.}
Our last example, due to K. Okada \cite{Ok},  concerns a  Hilbert modular newform $f$ on
$F=\Q(\sqrt{257})$ of weight $(2,2)$, level $1$, and Fourier coefficients in  $E=\Q(\sqrt{13})$.
Here the class number of $F$ is $3$.  

\begin{theo} \label{okada} Assume  $\ell\geq 5$ and $\ell\neq
  257$. Then the image of  $\overline{\rho}_{\lambda}$ 
is isomorphic to $\begin{cases}  \GL_2(\F_{\ell}) \text{, if } \ell\equiv \pm 1,\pm 3,\pm 4 \pmod{13}, \\
 \{\gamma\in \GL_2(\F_{\ell^2}),\det(\gamma)\in\F_{\ell}^{\times }\}\text{, otherwise}.\end{cases}$
\end{theo}

\noindent {\bf Proof : } Once again we have to discard  the reducible,  dihedral,
$A_4$, $S_4$, $A_5$  and inner twist cases.

$\relbar$ Assume that  $\overline{\rho}_{\lambda}$ is reducible :
$\overline{\rho}_{\lambda}^{\mathrm{s.s.}}=
\varphi\oplus\omega\varphi^{-1}$, where $\varphi$ is a mod $\lambda$ crystalline character of
$\G_F$. As $f$ is of level $1$, the same method as in  the previous section
shows that $\varphi$ is everywhere unramified (and therefore of  order dividing
the class number of $F$ - fact that we won't use).

Assume that $\ell\neq 61$. Let $\pi$ be one of the two prime ideals  of $F$ dividing $61$.
As  $\pi$ is  principal, 
by class field theory we deduce that $\varphi$ is trivial on  $\Frob_\pi$.
By evaluating $\overline{\rho}_{\lambda}^{\mathrm{s.s.}}= \varphi\oplus\omega\varphi^{-1}$
at $\Frob_\pi$ and taking traces, we obtain :
$$  -1\pm\sqrt{13} \equiv 1 + 61 = 62 \pmod{\lambda}.$$
 When  $\ell=61$, by applying the same argument with
a principal prime ideal  $\pi$ of $F$ dividing $67$, we get
 $$  -3\pm 3\sqrt{13} \equiv 1 + 67 = 68 \pmod{\lambda}.$$
Both these  congruences are impossible. Therefore
 $\overline{\rho}_{\lambda}$ is irreducible when $\ell\geq 5$ and $\ell\neq 257$.

\medskip
$\relbar$  As $f$ is of level $1$, and $F$ does not have unramified quadratic
extensions, it follows from Prop.\ref{image}(iii)   that the image of   $\overline{\rho}_{\lambda}$ cannot be of  dihedral  type for $\ell\geq 5$,  $\ell\neq 257$.

$\relbar$ The image of $\overline{\rho}_{\lambda}$ cannot be
of $A_4$, $S_4$ or $A_5$ type for  $\ell\geq 7$,  $\ell\neq 257$ 
(\cf  Prop.\ref{image}(ii)).
When $\ell=5$, we observe that the order of the image of
$\overline{\rho}_{5}(\Frob_{17A})$ in $\PGL_2(\F_{25})$ is at least $6$.

$\relbar$ As in the proof of Thm.\ref{cs} we see 
that $\overline{\rho}_{\lambda}$ cannot have inner twists. \hfill $\square$

\bigskip

\noindent Dept. d'Algebra i Geometria,
Universitat de Barcelona,
G.V. de les Corts Catalanes 585, 08007 Barcelona, Spain

\medskip
\noindent  e-mail : \texttt{ldieulefait@ub.edu}

\bigskip
\noindent UFR de Math{\'e}matiques, Universit{\'e} Paris 7,    Case 7012, 
2 place Jussieu,  75251 Paris, France

\medskip
\noindent  e-mail : \texttt{dimitrov@math.jussieu.fr}


\begin{thebibliography}{10}

\bibitem{carayol}
{\sc H.~Carayol}, {\em {Sur les repr{\'e}sentations $l$-adiques associ{\'e}es aux
  formes modulaires de Hilbert}}, Ann. Sci. {\'E}cole Norm. Sup., 19 (1986).

\bibitem{CoSc}
{\sc C.~Consani and J.~Scholten}, {\em {Arithmetic on a quintic threefold}},
  Internat. J. Math., 12 (2001), pp.~943--972.

\bibitem{Dem}
{\sc L.~Demb{\'e}l{\'e}}, {\em {Explicit computations of Hilbert modular forms on
  $\mathbb{Q}(\sqrt{5})$}}, PhD thesis, McGill University, Montreal, 2002.

\bibitem{Di}
{\sc L.~Dieulefait}, {\em {Newforms, inner twists, and the inverse Galois
  problem for projective linear groups}}, J. Th. Nombres Bordeaux, 13 (2001),
  pp.~395--411.

\bibitem{DV}
{\sc L.~Dieulefait and N.~Vila}, {\em {Projective linear gropus as Galois
  groups over $\mathbb{Q}$ via modular representations}}, J. Symb. Comp., 30
  (2000), pp.~799--810.

\bibitem{dim-hmv}
{\sc M.~Dimitrov}, {\em {Galois representations modulo $p$ and cohomology of
  Hilbert modular varieties}}.
\newblock preprint, 2004.

\bibitem{dim-ihara}
\leavevmode\vrule height 2pt depth -1.6pt width 23pt, {\em {On Ihara's lemma
  for Hilbert Modular Varieties}}.
\newblock preprint, 2004.

\bibitem{Ok}
{\sc K.~Okada}, {\em {Hecke Eigenvalues for Real Quadratic Fields}},
  Experimental Mathematics, 11 (2002), pp.~407--426.

\bibitem{ribet3}
{\sc K.~Ribet}, {\em {On $l$-adic representations attached to modular forms}},
  Invent. Math., 28 (1975), pp.~245--275.

\bibitem{ribet6}
\leavevmode\vrule height 2pt depth -1.6pt width 23pt, {\em {On $\ell$-adic
  representations attached to modular forms II}}, Glasgow Math. J., 27 (1985),
  pp.~185--194.

\bibitem{serre2}
{\sc J.-P. Serre}, {\em {Propri{\'e}t{\'e}s galoisiennes des points d'ordre fini des
  courbes elliptiques}}, Invent. Math., 15 (1972), pp.~259--331.

\bibitem{Ta}
{\sc R.~Taylor}, {\em {On Galois representations associated to Hilbert modular
  forms}}, Invent. Math., 98 (1989), pp.~265--280.

\bibitem{Ta3}
\leavevmode\vrule height 2pt depth -1.6pt width 23pt, {\em {On the meromorphic
  continuation of degree two $L$-functions}}.
\newblock preprint, 2001.

\end{thebibliography}
\end{document}